\documentclass[a4paper,11pt]{article}

\usepackage{amsmath,amsthm,amssymb,indentfirst}

\newcommand{\takes}{:}
\newcommand{\N}{\mathcal{N}}
\newcommand{\C}{\mathbb{C}}
\newcommand{\alg}{\mathcal{A}}

\newcommand{\hypalg}{\alg_h}
\newcommand{\fbalg}{\alg_p}
\newcommand{\FBalg}{\alg_{FB}}

\newcommand{\HS}{\mathcal{C}_2}
\newcommand{\R}{\mathbb{R}}

\newcommand{\halfplane}{\mathbb{H}}
\newcommand{\UHP}{\halfplane^+}
\newcommand{\LHP}{\halfplane^-}
\newcommand{\Q}{\mathbb{Q}}
\newcommand{\1}{\chi}
\newcommand{\bddop}{\mathcal{L}}
\newcommand{\calL}{\mathcal{L}}
\newcommand{\intint}{\int\!\!\!\int}
\DeclareMathOperator{\Int}{Int}
\DeclareMathOperator{\Alg}{Alg}
\DeclareMathOperator{\Lat}{Lat}
\DeclareMathOperator{\supp}{supp}
\DeclareMathOperator{\SOTcl}{SOT-cl}
\DeclareMathOperator{\Log}{Log}
\DeclareMathOperator{\Arg}{Arg}
\DeclareMathOperator{\sgn}{sgn}
\DeclareMathOperator{\clspan}{cl-span}

\def\swap#1#2{\let\tmp#1\let#1#2\let#2\tmp\let\tmp\relax}
\swap\phi\varphi

\makeatletter
\newcommand{\barenote}[1]{%
\renewcommand\@makefnmark{\relax}%
\addtocounter{footnote}{-1}%
\footnote{#1}%
\renewcommand\@makefnmark{\hbox{\@textsuperscript
{\normalfont\@thefnmark}}}%
}
\makeatother

\newtheorem{theorem}{Theorem}
\newtheorem{proposition}[theorem]{Proposition}
\newtheorem{lemma}[theorem]{Lemma}
\newtheorem{corollary}[theorem]{Corollary}

\theoremstyle{remark}
\newtheorem*{remark}{Remark}

\begin{document}
\title{Reflexivity of the translation-dilation algebras on $L^2(\R)$}
\author{R.~H.~Levene and
  S.~C.~Power} \maketitle

\begin{abstract}
  \barenote{2000 {\noindent\itshape Mathematics Subject
      Classification}. Primary 46K50.}\barenote{Keywords: reflexive algebra,
    translation invariant, Lie semigroup, hyperbolic algebra.}
  \barenote{The first author is supported by an EPSRC grant.}

  The hyperbolic algebra $\hypalg$, studied recently by Katavolos and
  Power~\cite{pow-kat:hyp}, is the weak star closed operator algebra
  on $L^2(\R)$ generated by $H^\infty(\R)$, as multiplication
  operators, and by the dilation operators $V_t$, $t\geq0$, given by
  $V_tf(x)=e^{t/2}f(e^t x)$. We show that $\hypalg$ is a reflexive
  operator algebra and that the four dimensional manifold
  $\Lat\hypalg$ (with the natural topology) is the reflexive hull of a
  natural two dimensional subspace.
\end{abstract}

\section{Introduction}
We resolve some problems posed in \cite{pow-kat:hyp} for the following
nonselfadjoint operator algebra $\hypalg$.  Let $M_\lambda$ and $V_t$
be the multiplication and dilation operators on $L^2(\R)$ given by
\[ M_\lambda f(x)=e^{i\lambda x}f(x),\qquad 
V_t f(x) = e^{t/2} f(e^t x)\]
and let $\hypalg$, the so-called hyperbolic algebra, be the weak-star
closed operator algebra generated by $\{M_\lambda, V_t\mid \lambda,
t\geq 0\}$. We can regard $\hypalg$ as a doubly nonselfadjoint algebra
in the sense that it is generated by the two copies of $H^\infty(\R)$
generated by these unitary semigroups separately. Under conjugation by
the Fourier-Plancherel transform the generators correspond to the
translations $D_\lambda f(x)=f(x-\lambda)$ and dilations (for $t\leq
0$) and for this reason we refer to $\hypalg$ as a
translation-dilation algebra. The weak star closed
translation-multiplication algebra, $\fbalg$ say, generated by $\{
M_\lambda, D_\mu\mid \lambda, \mu \ge 0\}$, is also doubly nonselfadjoint
and was analysed in~\cite{pow-kat:FB} and~\cite{power98:repns}. The
algebras $\fbalg$ and $\hypalg$ have similar properties and may be
viewed as fundamental examples of \emph{Lie semigroup algebras}. By
this we mean that they are weak star closed operator algebras
generated by the image of a Lie semigroup under a unitary
representation of the ambient Lie group. Indeed $\hypalg$ arises in
this way from the unitary representation of the $ax+b$ group given by
\[ \begin{pmatrix} a&b\\
  0&1\end{pmatrix} \mapsto M_bV_{\log a},\qquad a\neq 0,\ b\in \R,\] 
together with the Lie semigroup for $a\geq 1$, $b\geq0$.

It is a fundamental result of Sarason \cite{sarason72:weak} that
$H^\infty(\R)$ acting as a multiplication algebra on $H^2(\R)$ or
$L^2(\R)$ is a reflexive operator algebra. This means that the
invariant subspace lattice $\calL=\Lat \alg$ of the operator algebra
$\alg$ is such that
\[ \alg=\Alg\calL=\{A\mid AK\subseteq K\text{ for all $K\in \calL$}\}.\]
Since $\Alg\calL$ is the maximal operator algebra with lattice $\calL$ (and
hence reflexive) it is of particular interest to determine
reflexivity for operator algebras specified by generators, and indeed,
such a determination can be related to subtle problems of synthesis.
This is evident in Arveson's delicate analysis of operator algebra
synthesis and the density of pseudo-integral operators in CSL algebras
(see \cite{arveson74:invariant} and also
Froelich~\cite{froelich88:compact}). See also the slice map analysis
of the tensor product formula (and its failure) for various reflexive
algebras in \cite{kraus91:slice}, the analysis of reflexivity
for subnormal operators in \cite{thomson88:invariant}, and the
counterexamples of Wogen and Larson in \cite{wogen87:counterexamples}
and \cite{larson-wogen90:T+0}.

In \cite{pow-kat:FB} A.~Katavolos and the second author showed that
the translation-multiplication algebra on $L^2(\R)$ is reflexive by
identifying it with the (evidently) reflexive algebra
\[ \FBalg = \Alg(\N_a\cup \N_v) = (\Alg\N_a)\cap (\Alg\N_v)\]
where $\N_a=\{e^{i\lambda x}H^2(\R)\mid \lambda\in \R\}$ and $\N_v=
\{L^2[t,\infty]\mid -\infty\leq t\leq +\infty\}$ are the so-called
analytic nest and the Volterra nest respectively. The lattice
$\N_a\cup \N_v$, which we call the Fourier binest, is homeomorphic to
the unit circle when endowed with the natural topology, the relative
weak operator topology on the lattice of associated orthogonal
projections. Whilst this lattice provides the obvious subspaces in
$\calL_p=\Lat \fbalg$, the lattice $\calL_p$ turns out to be bigger and is
homeomorphic to the unit disc, with the Fourier binest as topological
boundary.  That is, $\calL_p$ is the reflexive hull of $\N_a\cup \N_v$;
\[\calL_p=\Lat \Alg(\N_a\cup \N_v).\]

The algebras $\fbalg$ and $\hypalg$ share several features over and
above being doubly nonselfadjoint. They are antisymmetric ($\alg\cap
\alg^*=\C I$) and they contain no non-trivial finite rank operators.
Both these facts rule out many of the basic techniques available for
the study of nest algebras (see~\cite{davidson:nest}). On the other
hand, $\fbalg$ and $\hypalg$ are not unitarily equivalent (nor
isometrically isomorphic) and $\calL_h=\Lat \hypalg$ is a 4-dimensional
Euclidean manifold which we indicate in Section~\ref{sec:reflxivity}.
The main results of the present paper show that $\hypalg$ is a
reflexive operator algebra and, in parallel with the comments above
for $\fbalg$, that $\calL_h$ is the reflexive hull of the sublattice
\[ \calL_M\cup \calL_S\]
where
\[\calL_M=\{L^2[-a,b]\mid a,b\in [0,\infty]\}\]
and
\[\calL_S = \{ |x|^{is}H^2(\R)\mid s\in \R\}\cup\{0,L^2(\R)\}.\]
Furthermore we obtain an explicit characterisation of the
Hilbert-Schmidt operators in $\hypalg$, and in fact this is the key to
the proof of reflexivity.

\section{Hilbert-Schmidt operators in $\Alg(\calL_M\cup\calL_S)$}

We write $\HS$ for the set of Hilbert-Schmidt operators on $L^2(\R)$,
so that $\HS=\{\Int k\mid k\in L^2(\R^2)\}$ where
\[(\Int k)f(x)=\int k(x,y)f(y)\,dy.\]

The following result gives a necessary support condition for $\Int k$
to be in $\Alg \calL_M$, and is simply proven. It is used several times
without reference when changing variables in what follows.
\begin{lemma}
  \label{lem:support}
  Let $\Int k\in \HS\cap \Alg\calL_M$. Then 
\begin{equation}
  \supp k\subseteq\{(x,y)\in
  \R^2\mid xy\geq0\text{ and }|y|\geq|x|\}.\tag*{\qedsymbol}
\end{equation}
\end{lemma}

The lemma shows that a Hilbert-Schmidt operator which leaves invariant
all the subspaces $L^2[-a,b]$ of $\calL_M$ has kernel function $k(x,y)$
supported in the cone bounded by the $y$-axis and the line $y=x$. We
now investigate the further constraints on $k(x,y)$ in order that
$\Int k$ should also leave invariant all the subspaces
$|x|^{is}H^2(\R)$ of the set $\calL_S$. We shall show that the condition
corresponds to the radially restricted functions $x\mapsto k(x,e^tx)$
lying in an appropriate function space for almost every $t$.  We begin
by defining this space, $V$.

Let $\UHP$ denote the open upper half-plane in $\C$ and $\LHP$ the
open lower half-plane. Let $\LHP_\Q$ be the set of points in $\LHP$
with rational real and complex parts. For $v,w\in \LHP$, set
\[ \tilde r_{v,w}(z)=\frac{1}{(z-v)(z-w)}. \]
Then $\tilde r_{v,w}$ is analytic in $\UHP$; let $r_{v,w}$ be the
boundary value function of $\tilde r_{v,w}$ on the real line.

The density assertion of the next lemma follows routinely from
Cauchy's theorem and the definition of $H^2(\R)$. It will be useful
for us in the proof of Proposition~\ref{prop:tensorsinhypalg}.

\begin{lemma}
  \label{lem:density}
  The set $\mathcal{R}_\Q=\{r_{v,w}\mid v,w\in \LHP_\Q\}$ has dense
  linear span in $H^2(\R)$.\qed
\end{lemma}

For $\alpha\in
\R$ and $z\in \UHP\cup \R\setminus\{0\}$, let
\[ \tilde p_\alpha(z)=\exp (\alpha\Log z)\]
where $\Log z = \log |z|+ i \Arg z $ and $\Arg z$ is the argument
function taking values in $[0,\pi]$.  Then $\tilde p_\alpha$ is
analytic in $\UHP$. Let $p_\alpha$ be the restriction of $\tilde
p_\alpha$ to $\R\setminus\{0\}$ and let $p=p_{1/2}$ and $q=p_{-1/2}$.
Note that $\overline{p(x)}=\sgn(x)p(x)$ and
$\overline{q(x)}=\sgn(x)q(x)$ where $\sgn(x)=x/|x|$ for $x\neq 0$.

Now define the space $V$ to be the closed linear span in
$L^2(|x|\,dx)$ of the set $\{ \overline{q} r_{v,w}\mid v,w\in
\LHP_\Q\}$. Observe that the unitary $M_{\overline{p}}\takes
L^2(|x|\,dx)\to L^2(\R)$ restricts to a unitary map $V\to
H^2(\R)$, and $M_{\overline{p}}^*=M_{\overline{q}}$. This shows that,
writing $\clspan{A}$ for the closed linear span of $A$, we have
\begin{align*}
  V^\perp&=
  M_{\overline{p}}^*\big(M_{\overline{p}}V\big)^\perp \\
  &= M_{\overline{q}}\big(H^2(\R)\big)^\perp\\
  &= M_{\overline{q}}\overline{H^2(\R)}\\
  &= M_{\overline{q}} \clspan{\{\overline{r_{v,w}}\mid v,w\in \LHP_\Q\}}\\
  &= M_{\sgn}\clspan{\{ q\overline{r_{v,w}}\mid  v,w\in \LHP_\Q\}}\\
  &= M_{\sgn}\overline{V}.
\end{align*}

Henceforth let $W=L^2(\R_+,e^t\,dt)$. We will
write $\|\cdot\|_V$ and $\|\cdot\|_W$ for norms taken in the spaces
$L^2(|x|\,dx)$ and $L^2(e^t\,dt)$ respectively.

The next lemma is elementary.

\begin{lemma}
  \label{lem:lebesgue}
  Let $\mu$ and $\nu$ be measures on $\R$, and let $M$ and $N$ be
  closed subspaces of $L^2(\mu)$ and $L^2(\nu)$ respectively. Let
  $f\in L^2(\mu\times \nu)$, and suppose that for almost every $y$,
  $x\mapsto f(x,y)\in M$, and for almost every $x$, $y\mapsto
  f(x,y)\in N$.  Then $f\in M\otimes N$, the tensor product Hilbert
  space.\qed
\end{lemma}

We now obtain a necessary condition on the kernel functions $k(x,y)$
of operators in $\Alg(\calL_M\cup \calL_S)$. Here and elsewhere we make the
natural identification between the tensor product $M\otimes N$ and the
closed linear span of $\{ m(x)n(y)\mid m\in M,\ n\in N\}$ in
$L^2(\mu\times \nu)$.

\begin{proposition}
  \label{prop:kernelsinTP}
  Let $\Int k$ be a Hilbert-Schmidt operator leaving invariant each of
  the subspaces $L^2[-a,b]$, $a,b\geq0$ and $|x|^{is}H^2(\R)$, $s\in
  \R$. Then the function $(x,t)\mapsto k(x,e^tx)$ is in the tensor
  product Hilbert space $V\otimes W$.
\end{proposition}
\begin{proof}
  We aim to apply Lemma~\ref{lem:lebesgue}. Firstly, for almost every
  $x$, $t\mapsto k(x,e^t x)\in W$, since
  \[\intint |k(x,y)|^2\,dx\,dy = 
  \intint |k(x,e^tx)|^2e^t|x|\,dx\,dt<\infty.\] This also shows that for
  almost every $t$, $x\mapsto k(x,e^tx)$ is in $L^2(|x|\,dx)$. It
  remains to show that this function lies in $V$.
  
  Let $X=\Int k$.  Then for every $s\in \R$ we have
  $X|x|^{is}H^2(\R)\subseteq |x|^{is}H^2(\R)$, so for every $h_1$ and
  $h_2$ in $H^2(\R)$,
  \begin{align*}
    0&=\langle X|x|^{is}h_1,|x|^{is}\overline{h_2}\rangle\\
    &=\int\left[\int k(x,y)|y|^{is}h_1(y)\,dy\right]\,|x|^{-is}
    h_2(x)\,dx.
  \end{align*}
  Now using Lemma~\ref{lem:support} we make the change of variables
  $y=e^tx$ to get
  \begin{equation}
    0 =\intint e^{t/2}|x|k(x,e^tx)
    e^{ist}(V_th_1)(x)h_2(x)\,dt\,dx.\label{eq:T1}
  \end{equation}
  Now
  \begin{align}\intint |e^{t/2}xk(x,e^tx)(h_2V_th_1)(x)|\,dx\,dt 
    &=\intint |k(x,y)h_2(x)h_1(y)|\,dx\,dy\notag\\ &= \langle
    (\Int|k|)|h_1|,|h_2|\rangle <\infty, \label{eq:fubini}
  \end{align}%
  so the integrand in (\ref{eq:T1}) is in $L^1(dt\times dx)$ and we
  can apply Fubini's theorem to get
  \begin{equation}
    0=\int \left\{ \int e^{t/2}|x|k(x,e^tx)(h_2V_th_1)(x)\,dx\right\}
    \,e^{ist}\,dt.  \label{eq:T}
  \end{equation}
  Let $T(t)$ be the expression in curly brackets. By virtue of
  (\ref{eq:fubini}) $T$ is in $L^1(\R)$ and so \eqref{eq:T} can be
  interpreted as saying that the Fourier transform of the $L^1$
  function $T$ is zero. Hence for almost every $t$, $T(t)=0$.
  
  For fixed $t$, as $h_1$ and $h_2$ run over $H^2(\R)$, $h_2V_th_1$
  runs over $H^1(\R)$, so we can pick $h_1$ and $h_2$ such that
  $h_2V_th_1$ runs over a countable dense subset of $H^2(\R)$.  Since
  $T=0$, for almost every $t$ we have
  \[\int k(x,e^tx) f(x)\,|x|dx=0\]
  for every $f\in H^2(\R)$. We can now approximate again: since
  $qr_{v,w}$ is in $H^2(\R)$ for each $v,w\in \LHP$, for almost every
  $t$ we have
  \[ 0=\int k(x,e^tx) q(x)r_{v,w}(x) \,|x|\,dx
  = \langle x\mapsto k(x,e^tx), \overline{qr_{v,w}}\rangle_{V}
  \]
  for each $v,w\in \LHP_\Q$. But $\{\overline{qr_{v,w}}\mid v,w\in
  \LHP_\Q\}$ is countable and has $\|\cdot\|_V$-dense linear span in
  $V^\perp$, so we conclude that for almost every $t$, $x\mapsto
  k(x,e^tx)$ is in $V$.
    
  Lemma~\ref{lem:support} and a simple computation show that for
  almost every $x$, $t\mapsto k(x,e^tx)$ is in $W$.  Now
  Lemma~\ref{lem:lebesgue} ensures that $k(x,e^tx)$ is in $V\otimes
  W$.
\end{proof}

\section{Hilbert-Schmidt operators in $\hypalg$}
\label{sec:HS}

We have the elementary inclusions
\begin{equation} 
  \label{eq:inclusions}
  \hypalg\subseteq \Alg(\Lat\hypalg)\subseteq \Alg(\calL_M\cup\calL_S).
\end{equation}
To show that we in fact have equalities, and in particular that
$\hypalg$ is reflexive, we shall establish two facts. Firstly, that
the necessary condition of Proposition~\ref{prop:kernelsinTP} is in
fact sufficient for $\Int k$ to belong to $\hypalg$. Thus, $\hypalg$
and $\Alg(\calL_M\cup \calL_S)$ will contain the same Hilbert-Schmidt
operators.  Secondly we show that for each algebra the Hilbert-Schmidt
operators are dense for the weak star topology.

\begin{lemma}
  \label{lem:Vphi}
  Let $\phi\in W$ have compact support. Then for $f\in L^2(\R)$ with
  compact support and almost every $x\in \R$, the integral
  \[ \tilde V_\phi f(x)=\int e^{t/2}\phi(t)V_tf(x)\,dt,\]
  converges, and $\tilde V_\phi$ extends to a bounded operator
  $V_\phi\in \bddop(L^2(\R))$. Moreover, $V_\phi\in \hypalg$.
\end{lemma}

\begin{proof}
  Checking that $\tilde V_\phi f(x)$ converges for almost every $x$
  when $f$ is compactly supported is easily done. Let $f$ and $g$ be
  in $L^2(\R)$. Then
  \begin{align*}
    \bigg|\intint e^{t/2}\phi(t)V_tf(x)\overline{g(x)}\,dx\,dt\bigg|
    &\leq \int |e^{t/2}\phi(t) \langle V_tf,g\rangle|\,dt
    \\ &\leq \|f\|\,\|g\|\,\|\phi\|_{L^1(e^{t/2}\,dt)}.
  \end{align*}
  Since $\phi$ is compactly supported and in $W$, it is in
  $L^1(e^{t/2}\,dt)$ and so the sesquilinear form
  \[ \langle f,g\rangle_\phi = \intint
  e^{t/2}\phi(t)V_tf(x)\overline{g(x)}\,dx\,dt \] is bounded and defines
  a bounded operator $V_\phi$ extending $\tilde V_\phi$.
  
  To show that $V_\phi\in \hypalg$, without loss of generality we can
  assume that the support of $\phi$ is contained in the interval
  $[0,1]$. For positive integers $m<n$, let
  \[\alpha_{n,m}=\int_{m/n}^{(m+1)/n} e^{s/2}\phi(s)\,ds,\]
  and let $T_n=\sum_{m=0}^{n-1} \alpha_{n,m}V_{m/n}$. Then $T_n\in
  \hypalg$ for every $n$. We claim that $T_n\to V_\phi$ in WOT. For if
  $f,g\in L^2(\R)$, then
  \begin{align*}
    \langle (V_\phi-T_n)f,g\rangle 
    &=\int \overline{g(x)} \biggl(\int
    e^{t/2}\phi(t)(V_tf)(x)\,dt-\sum_{m=0}^{n-1}\alpha_{n,m}(V_{m/n}f)(x)\biggr)\,dx\\
    &=\int \overline{g(x)} \sum_{m=0}^{n-1}\int_{m/n}^{(m+1)/n}
    e^{t/2}\phi(t)\big( (V_tf)(x)-(V_{m/n}f)(x)\big)\,dt\,dx\\
    &=\int e^{t/2}\phi(t)\big(\langle V_tf,g\rangle-\langle
    V_{m(n,t)/n} f,g\rangle\big)\,dt
  \end{align*}
  where $m(n,t)=\lfloor nt\rfloor$. Now for any $t$, $|\langle
  V_tf,g\rangle|\leq \|f\|\,\|g\|$, and $e^{t/2}\phi(t)$, being a
  compactly supported element of $L^2(\R)$, is in $L^1(\R)$. Moreover,
  $V_{m(n,t)/n} \to V_t$ in WOT as $n\to \infty$. So by dominated
  convergence, $\langle (V_\phi-T_n)f,g\rangle\to 0$ as $n\to \infty$.
  For every $n$, $\|T_n\|\leq \|\phi\|_{L^1(e^t\,dt)}$, so $(T_n)$ is
  a bounded sequence in $\bddop(L^2(\R))$. The
  $\mathrm{w}^*$-topology and WOT agree on bounded sets; since
  $\hypalg$ is $\mathrm{w}^*$-closed, this shows that $V_\phi\in
  \hypalg$.
\end{proof}

Given a function $F\takes \R\times \R_+\to \C$, we define the function
$k_{F}\takes \R^2\to \C$ by
\[   k_{F}(x,y)=
  \begin{cases}
    F(x,\log (y/x))& \text{if $xy>0$}\\
    0 & \text{otherwise}
  \end{cases}.
\]
We also write $k_{h,\phi}=k_{h(x)\phi(t)}$ when $h\takes \R\to \C$ and
$\phi\takes\R_+\to \C$. Note that this means that $k_{h,\phi}(x,e^tx)
= h(x)\phi(t)$ for $t>0$.

\begin{proposition}\label{prop:tensorsinhypalg}
  Let $h\in V$ and $\phi\in W$. Then $\Int k_{h,\phi} \in \hypalg\cap
  \HS$.
\end{proposition}
\begin{proof}
  Let $k=k_{h,\phi}$. Since
  \begin{align}
    \|k\|_{L^2(\R^2)}^2
    &=\intint |k(x,y)|^2\,dx\,dy\notag\\
    &=\int |\phi(t)|^2e^t\,dt\cdot\int |h(x)|^2|x|\,dx\notag\\
    &=\|\phi\|_{W}^2\|h\|_{V}^2,
    \label{eq:normVW}
  \end{align}
  we have $k\in L^2(\R^2)$ and so $\Int k\in \HS$. Since
  $\|k\|_{L^2(\R^2)}=\|\Int k\|_{\HS}$, (\ref{eq:normVW}) also shows
  that if $h_n\to h$ in $\|\cdot\|_V$, then the operators $\Int
  k_{h_n,\phi}$ converge to $\Int k$ in Hilbert-Schmidt norm, and so
  in operator norm. Since $\hypalg$ is closed in operator norm, by the
  definition of $V$ we may assume that $h=\overline{q}r_{v,w}$ for
  some $v,w\in \LHP_\Q$. Also, if $\phi$ has compact support, then
  \begin{align*}
    (\Int k)f(x) &= \int_{\{y\mid xy>0\}}
    h(x)\phi(\log(y/x))f(y)\,dy\\
    &=|x|h(x)\int e^{t/2} \phi(t)V_tf(x)\,dt\\
    &= M_gV_\phi f(x),
  \end{align*}
  where $g(x)=|x|h(x) = p(x)r_{v,w}(x)\in H^\infty(\R)$, as is easily verified.
  So $\Int k=M_g V_\phi$. Now $M_g\in \hypalg$ since $g\in
  H^\infty(\R)$. By Lemma~\ref{lem:Vphi}, $V_\phi\in \hypalg$ too,
  so their product $\Int k$ is also in $\hypalg$.  
  
  If $\phi$ is not compactly supported, we can simply restrict $\phi$
  to compact subsets of $\R_+$ to get a sequence of Hilbert-Schmidt
  operators in $\hypalg$ that converges in Hilbert-Schmidt norm to
  $\Int k$ by (\ref{eq:normVW}) again, and so in operator norm. So
  $\Int k$ is in $\hypalg\cap \HS$ as desired.
\end{proof}

\begin{corollary}\label{cor:HS-hypalg}
  $\hypalg\cap \HS = \Alg(\calL_M\cup \calL_S)\cap \HS = \{ \Int k_{F}\mid
  F\in V\otimes W \}$.
\end{corollary}
\begin{proof}
  Let $F\in V\otimes W$; then $F=\sum \lambda_i h_i\otimes \phi_i$
  with convergence in $L^2(|x|\,dx\times e^t\,dt)$. So $\Int
  k_{F}=\sum \lambda_i \Int k_{h_i, \phi_i}$ with convergence in
  Hilbert-Schmidt norm. By Proposition~\ref{prop:tensorsinhypalg},
  each $\Int k_{h_i, \phi_i}$ is in $\hypalg\cap \HS$, so $\Int
  k_{F}$ is too.
  
  Conversely, if $\Int k\in \Alg (\calL_M\cup\calL_S)\cap \HS\supseteq
   \hypalg\cap \HS$ then Proposition~\ref{prop:kernelsinTP} implies
  that $k=k_{F}$ for some $F\in V\otimes W$ and the corollary follows.
\end{proof}

\begin{proposition}\label{prop:BAI}
  $\hypalg\cap \HS$ contains a bounded approximate identity---that is,
  a sequence that is bounded in operator norm and converges in the
  strong operator topology to the identity.
\end{proposition}
\begin{proof}
  Let
  \[h_n(x)=\frac{-n^2\sgn(x)}{(x+n^{-1}i)(x+ni)^2},\qquad
  \phi_n(t)=ne^{-t/2}\1_{[0,n^{-1}]}(t)\] %
  and let $g_n(x)=|x|h_n(x)$.  Then $|g_n(x)|\leq 1$ for each real $x$
  and $g_n(x)\to 1$ uniformly on compact sets, so $M_{g_n}\to I$ in
  SOT.  Since $\overline{p}h_n\in H^2(\R)$, we have $h_n\in V$.
  Moreover, $\phi_n\in W$, and a calculation using Jensen's inequality
  shows that $\|V_{\phi_n}\|\leq 1$ for each $n$. If we take $f\in
  L^2(\R)$ continuous with compact support then by dominated
  convergence,
  \[ \|V_{\phi_n}f-f\|^2 =\int\Big|\int_{0}^{1/n} nV_tf(x)\,dt -
  f(x)\Big|^2\,dx\to 0 \text{ as }n\to \infty.\] Since such functions
  are dense in $L^2(\R)$, $V_{\phi_n}\to I$ boundedly in SOT. Let
  $k_n=k_{h_n,\phi_n}$. Then, exactly as in the proof of
  Proposition~\ref{prop:tensorsinhypalg}, $\Int
  k_n=M_{g_n}V_{\phi_n}\in \hypalg\cap \HS$ and so $\Int k_n\to I$
  in SOT.
\end{proof}

\begin{remark}
  The bounded approximate identity above is considerably simpler than
  that obtained in \cite{pow-kat:hyp}.
\end{remark}\medskip

\section{Reflexivity}
\label{sec:reflxivity}

In \cite{pow-kat:hyp}, A.~Katavolos and the second author showed that
the lattice of invariant subspaces for $\hypalg$ is 
\[ \calL_h = \Lat \hypalg = 
\Bigg( \bigcup_{s\in \R,\,|\theta|=1}\calL_{s,\theta}\Bigg)
\cup \calL_M\] where
\begin{align*}
  \calL_{s,\theta}&=\{u_\theta(x)|x|^{is}e^{i(\lambda x+\mu
    x^{-1})}H^2(\R)\mid \lambda,\mu\geq 0\}.
\end{align*}
Here $u_\theta(x)$ is the function taking the value $1$ for $x>0$ and
$\theta$ for $x<0$. Also it was shown that with the natural topology
$\calL_h$ is connected and is a $4$-dimensional Euclidean manifold.

\begin{theorem}
  The hyperbolic algebra $\hypalg$ is reflexive.
\end{theorem}
\begin{proof}
  Take the sequence $X_n$ in $\hypalg\cap \HS$ from
  Proposition~\ref{prop:BAI} such that $\|X_n\|\leq 1$ for every $n$
  and $X_n\to I$ strongly. If $T$ is an operator in $\hypalg$, then
  $T$ is the SOT limit of the operators $TX_n\in \hypalg\cap \HS$.
  Similarly, if $S$ is an operator in $\Alg\calL_h$, then $S$ is the SOT
  limit of the operators $SX_n\in \Alg(\calL_M\cup\calL_S)\cap \HS$.  Since
  $\hypalg$ and $\Alg(\calL_M\cup\calL_S)$ are both strongly closed and by
  the results of Section~\ref{sec:HS} they contain the same
  Hilbert-Schmidt operators, upon taking strong operator topology
  closures, we see that
  \[\hypalg=\SOTcl(\hypalg\cap \HS)=\SOTcl\big(\Alg(\calL_M\cup\calL_S)\cap \HS\big)
  =\Alg(\calL_M\cup\calL_S).\] So the inclusions in \eqref{eq:inclusions}
  are equalities and $\hypalg$ is reflexive.
\end{proof}

We now obtain some consequences for sublattices of $\calL_h$. Note that
part~(ii) of the next theorem is analogous to the fact that the
lattice $\calL_p$ of the translation-multiplication algebra $\FBalg$ is
the reflexive hull of the Fourier binest $\N_a\cup \N_v$.

\begin{theorem}
  Let $\calL_M$, $\calL_V$ and $\calL_S$ be the subspace lattices
  \begin{align*}
    \calL_M&=\{L^2[-a,b]\mid a,b\in [0,\infty]\},\\
    \calL_V&=\{u_\theta(x)|x|^{is}H^2(\R)\mid \theta\in S^1,s\in \R\}\cup\{0,L^2(\R)\},\\
    \calL_S &= \{ |x|^{is}H^2(\R)\mid s\in \R\}\cup\{0,L^2(\R)\}.
  \end{align*}
  
  (i) The lattice $\calL_M$ (respectively $\calL_V$) is the set of subspaces
  of $\calL_h$ which are reducing subspaces for the semigroup
  $\{M_\lambda\mid \lambda\in \R_+\}$ (respectively $\{V_t\mid t\in
  \R_+\}$).
  
  (ii) $\calL_h$ is the reflexive hull of $\calL_M\cup \calL_V$.
  
  (iii) $\calL_h$ is the reflexive hull of $\calL_M\cup \calL_S$.
\end{theorem}
\begin{proof}
  Part~(i) is elementary. Part~(iii) is contained in the last proof
  and (ii) follows from~(iii).
\end{proof}

\newpage
\bibliographystyle{plain}

{\footnotesize
\bibliography{Ah-refl}

\noindent{\sc Department of Mathematics and Statistics, Lancaster University,
  LA1 4YF, England}

\noindent{\it Fax}: +44 1524 592681

\noindent{\it E-mail addresses}: {\tt r.levene@lancaster.ac.uk}, {\tt s.power@lancaster.ac.uk}
}
\end{document}